\RequirePackage{ifpdf}
\ifpdf 
\documentclass[pdftex]{sigma}
\else
\documentclass{sigma}
\fi

\begin{document}

\allowdisplaybreaks

\renewcommand{\PaperNumber}{104}

\FirstPageHeading

\ShortArticleName{Note on Dilogarithm Identities}

\ArticleName{Note on Dilogarithm Identities\\ from Nilpotent Double Af\/f\/ine Hecke Algebras}

\Author{Tomoki NAKANISHI}

\AuthorNameForHeading{T.~Nakanishi}

\Address{Graduate School of Mathematics, Nagoya University, Chikusa-ku, Nagoya, 464-8604, Japan}
\Email{\href{mailto:nakanisi@math.nagoya-u.ac.jp}{nakanisi@math.nagoya-u.ac.jp}}

\ArticleDates{Received November 15, 2012, in f\/inal form December 22, 2012; Published online December 25, 2012}

\Abstract{Recently Cherednik and Feigin [arXiv:1209.1978] obtained several Rogers--Ra\-ma\-nujan type identities
via the nilpotent double af\/f\/ine Hecke algebras (Nil-DAHA).
These identities further led to a series of dilogarithm identities,
some of which are known,
while some are left conjectural.
We conf\/irm and explain all of them
by showing the connection
with $Y$-systems
associated with (untwisted and twisted) quantum af\/f\/ine Kac--Moody algebras.}


\Classification{17B37; 13F60}

\section{Dilogarithm identities from Nil-DAHA}

Let $R_n$ be a root system of f\/inite type and of rank $n$
with non-degenerate bilinear form $(\ , \ )$,
and let $\alpha_i$ and $\omega_i$ be the simple roots and
the fundamental weights of $R_n$.
The Cartan matrix $C=(c_{ij})_{i,j=1}^n$
is given by
$c_{ij}=2(\alpha_i,\alpha_j)/(\alpha_i,\alpha_i)$.
Following \cite{Cherednik12},
let $A=(a_{ij})$ and
$A^{\flat}=(a^{\flat}_{ij})$ be the matrices with
  $a_{ij}=2(\omega_i,\omega_j)$ and
  $a^{\flat}_{ij}=(\omega_i,\omega_j)$, respectively.
Set
$\nu_i = (\alpha_i,\alpha_i)/2$.
Then,
$\nu_i^{-1}(\alpha_i,\omega_j)={\delta}_{ij}$,
and we have
\begin{gather}
\label{eq:Ainv}
\big(A^{\flat}\big)^{-1}
 = \big(c_{ij}\nu_j^{-1}\big)_{i,j=1}^n.
\end{gather}
Below we normalize the bilinear form as
$(\alpha_{\mathrm{short}},\alpha_{\mathrm{short}})=2$
so that $\nu_i
 \in \{1,2,3\}$.

Let $L(x)$ be the Rogers dilogarithm function
\begin{gather*}
L(x)=
-\frac{1}{2}
\int_{0}^x
\left\{
\frac{\log (1-y)}{y}
+
\frac{\log y} {1-y}
\right\}
{\mathrm d}y.
\end{gather*}
In \cite[equation~(3.34)]{Cherednik12} Cherednik and Feigin presented
two (partially conjectural)
series of di\-loga\-rithm identities.
Let $A'=(a'_{ij})_{i,j=1}^n$ be either $A$ or $A^{\flat}$ as above.
Let $Q_i$ ($i=1,\dots,n$) be the unique solution of the system of equations
\begin{gather}
\label{eq:Qeq}
(1-Q_i)^{\nu_i} = \prod_{j=1}^n Q_j^{a'_{ij}}
\end{gather}
in the range $0< Q_i <1$.
Then, the following identity was proposed
\begin{gather}
\label{eq:CFdilog}
\frac{6}{\pi^2}
\sum_{i=1}^n \nu_i L(Q_i) = L_{A'},
\end{gather}
where the value $L_{A'}$ is the rational number given in
Table~\ref{tab:L}.
\begin{table}[t]\centering
\renewcommand\arraystretch{1.6}
\renewcommand\tabcolsep{5pt}
\caption{The value $L_{A'}$.}\label{tab:L}
\vspace{1mm}
\begin{tabular}{c|cccccccccc}
\hline
$R_n$ & $A_n$ &$B_n$ &$C_n$ &$D_n$ &$E_6$ &$E_7$ &$E_8$ &
$F_4$ &$G_2$ & $T_n$\\
\hline
$L_{A}$ & $\frac{n(n+1)}{n+3}$ &
$\frac{n(2n-1)}{n+1}$ &$n$ &$n-1$ &
$\frac{36}{7}$ &$\frac{63}{10}$ &
$\frac{15}{2}$ &
$\frac{36}{7}$ &$3$ & $\frac{n(2n+1)}{2n+3}$
\\
\hline
$L_{A^{\flat}}$ & $\frac{n(n+1)}{n+4}$
 &$\frac{2n(2n-1)}{2n+3}$
 &$\frac{2n(n+1)}{2n+3}$ &$\frac{2(n-1)n}{2n+1}$ &
$\frac{24}{5}$ &$6$ &$\frac{80}{11}$ &
$\frac{24}{5}$ &$\frac{8}{3}$
&$\frac{n(2n+1)}{2n+4}$\\
\hline
\end{tabular}
\end{table}

In addition, there are identities for `type $T_n$' (tadpole type).
We def\/ine the `Cartan mat\-rix'~$C$ as almost the same as
type $A_n$ except that the last diagonal entry is~1, not~2.
Also the matri\-ces~$A$ and~$A^{\flat}$
are def\/ined by $a_{ij}=2\min(i,j)$ and $a^{\flat}_{ij}=\min(i,j)$,
respectively; the latter is the same matrix $A^{\flat}$ for  type $C_n$.
Then, $(A^{\flat})^{-1}=C$ holds.
We
set $\nu_i=1$.
Again, \eqref{eq:CFdilog} should hold
for the value $L_{A'}$  in
Table \ref{tab:L}, where $T_n$ is formally included  as
a member of $R_n$.

In \cite{Cherednik12} these identities were partially obtained and generally
 motivated by the Rogers--Rama\-nujan type
identities arising from nilpotent af\/f\/ine Hecke algebras (Nil-DAHA),
but only some of them are identif\/ied with the known identities.

The authors
of  \cite{Cherednik12}
expected the connection between~\eqref{eq:CFdilog} and
dilogarithm identities from some $Y$-systems (and cluster algebras behind them).
In this note we answer this question af\/f\/irmatively,
and, in particular, we conf\/irm all the identities in question.
The note has considerable overlap with the paper by Lee~\cite{Lee11},
but it is written for a~dif\/ferent purpose and in a~dif\/ferent perspective.

\section[Dilogarithm identities for $Y$-systems of simply laced type]{Dilogarithm identities for $\boldsymbol{Y}$-systems of simply laced type}

Let us recall the following dilogarithm identities proved by cluster algebra
method~\cite{Nakanishi09}. For $\ell=2$ see also \cite{Chapoton05}.

Let $C$ be any Cartan matrix of simply laced type $R_n=A_n,\,D_n,\,E_6,\,E_7,\,E_8$, and let $\ell\geq 2$  be any integer (called the {\em level}).
Let $Y_m^{(i)}$ ($i=1,\dots,n$; $m=1,\dots, \ell-1$)
 be the unique real positive solution of the system of equations
\begin{gather}
\label{eq:Y-ADE}
\big(Y^{(i)}_m\big)^2=
\frac{
\prod\limits_{j=1}^n \big(1+Y_m^{(j)}\big)^{2\delta_{ij}-c_{ij}}}
{
\big(1+Y_{m-1}^{(i)}{}^{-1}\big)
\big(1+Y_{m+1}^{(i)}{}^{-1}\big)
},
\end{gather}
where $Y_{0}^{(i)}{}^{-1}=Y_{\ell}^{(i)}{}^{-1}=0$.
\begin{theorem}[{\cite[Corollay 1.9]{Nakanishi09}}]
\label{thm:cDI1}
The following identity holds
\begin{gather}
\label{eq:dilogADE}
\frac{6}{\pi^2}
\sum_{i=1}^n \sum_{m=1}^{\ell-1}L\left(\frac{Y^{(i)}_m}{1+Y^{(i)}_m}\right) =
\frac{
(\ell-1)nh
}
{h+\ell
},
\end{gather}
where $h$ is the Coxeter number of type $R_n$, i.e.,
$n+1$, $2n-2$, $12$, $18$, $30$ for
$A_n$, $D_n$, $E_6$,
$E_7$, $E_8$, respectively.
\end{theorem}

The system of equations~\eqref{eq:Y-ADE} is called the {\em
 level $\ell$ constant $Y$-system}
associated with the quantum af\/f\/ine Kac--Moody algebra
of (untwisted) type $R^{(1)}_n$,
which is a specialization of the cor\-responding (non-constant) $Y$-system.
It is also known as the (constant) $Y$-system of \mbox{$R_n \times A_{\ell-1}$}.
See~\cite{Inoue10c} for more information.

We explain below that all the identities \eqref{eq:CFdilog} in question
are the ones in \eqref{eq:dilogADE} for $\ell=2,3$, or their specializations.

\section[Identification with $Y$-systems from quantum affine Kac-Moody algebras]{Identif\/ication with $\boldsymbol{Y}$-systems\\ from quantum af\/f\/ine Kac--Moody algebras}

\subsection[The non-$\flat$ case]{The non-$\boldsymbol{\flat}$ case}

We consider the case $A'=A$.
We use the change of variables
\begin{gather}
\label{eq:QY}
Q_i=\frac{Y_i}{1+Y_i},
\end{gather}
so that the range $0<Q_i<1$ corresponds to the range $0< Y_i$.
Then, using \eqref{eq:Ainv}, one can transform the equations~\eqref{eq:Qeq} into the form
\begin{gather*}
\prod_{j=1}^n
\left(
 \frac{1}{1+Y_j}
\right)^{c_{ij}}
=
\left(
\frac{Y_i}
{1+Y_i}\right)^2,
\end{gather*}
which is equivalent to
\begin{gather}
\label{eq:Ynflat}
Y_i^2=
\prod_{j=1}^n (1+Y_j)^{2\delta_{ij}-c_{ij}}.
\end{gather}

For simply laced types $A_n$, $D_n$, $E_6$,
$E_7$, $E_8$, \eqref{eq:Ynflat}  coincides with
the {\em level}~2
constant $Y$-system~\eqref{eq:Y-ADE} of untwisted type~$R_n^{(1)}$
by identifying $Y_i$ with $Y^{(i)}_1$.
Then, the right hand side of~\eqref{eq:dilogADE} with $\ell=2$ gives the value of~$L_A$,
agreeing with Table~\ref{tab:L}.

In contrast, for types
$B_n$, $C_n$, $F_4$, $G_2$,
\eqref{eq:Ynflat} coincides with
the level~2 constant $Y$-system of
{\em twisted type $S^{(r)}_m$}
(in the sense of  \cite[Remark~9.22]{Inoue10c}),
where $S^{(r)}_m$ is the
Langlands dual of $R^{(1)}_n$.
See Table~\ref{tab:dual} for the Langlands dual of af\/f\/ine type.
Also see \cite[Section~9]{Inoue10c} for the full version of $Y$-systems of twisted type.
In this case, the direct inspection of the Cartan matrix shows that
the equation~\eqref{eq:Ynflat} can be obtained from
the level~$2$ constant $Y$-system of
the {\em untwisted} (and simply laced)
type~$S^{(1)}_m$ by {\em folding}, i.e., identifying the variables
with the diagram automor\-phism~$\sigma$ of~$S_m$.
This is possible, due to the symmetry $Y^{(i)}_m\leftrightarrow
Y^{(\sigma(i))}_m$ of the  $Y$-system~\eqref{eq:Y-ADE}.
Furthermore, it is easy to see that~$\nu_i$ coincides with
the number of elements in the $\sigma$-orbit of~$i$.
Thus, we obtain the identity~\eqref{eq:CFdilog} for type~$R_n$
with $L_A(R_n)= L_A(S_m)$.
For example, for $R_n=B_n$, $L_A(B_n)=L_A(A_{2n-1})$.
This conf\/irms and explains Table~\ref{tab:L}.

\begin{table}[t]\centering
\renewcommand\arraystretch{1.6}
\renewcommand\tabcolsep{5pt}
\caption{The Langlands dual of af\/f\/ine type.} \label{tab:dual}
\vspace{1mm}

\begin{tabular}{c|ccccccccc}
\hline
$R_n^{(1)}$ & $A_n^{(1)}$ &$B_n^{(1)}$ &
$C_n^{(1)}$ &$D_n^{(1)}$ &$E_6^{(1)}$ &
$E_7^{(1)}$ &$E_8^{(1)}$ &
$F_4^{(1)}$ &$G_2^{(1)}$\\
\hline
$S_m^{(r)}$ & $A_n^{(1)}$ &$A_{2n-1}^{(2)}$ &
$D_{n+1}^{(2)}$ &$D_n^{(1)}$ &$E_6^{(1)}$ &
$E_7^{(1)}$ &$E_8^{(1)}$ &
$E_6^{(2)}$ &$D_4^{(3)}$\\
\hline
\end{tabular}
\end{table}

Finally, for type $T_n$,
\eqref{eq:Ynflat} coincides with
the level 2 constant $Y$-system
of type $A^{(2)}_{2n}$
(in the sense of
\cite[Remark 9.22]{Inoue10c}).
Note that $A^{(2)}_{2n}$ is self-dual under the Langlands duality.
Again, this $Y$-system is obtained by the folding of level $2$ constant $Y$-system
of {\em untwisted} $A^{(1)}_{2n}$. Since we set $\nu_i=1$, the multiplicities
are discarded in~\eqref{eq:CFdilog}.
Therefore, $L_A(T_n)=L_A(A_{2n})/2$.
Actually, this connection is known in~\cite{Cherednik12} and other literature.

\subsection[The $\flat$ case]{The $\boldsymbol{\flat}$ case}

We consider the case $A'=A^{\flat}$.
By the same change of variables~\eqref{eq:QY},
one can transform the equations~\eqref{eq:Qeq} into the form
\begin{gather*}
\prod_{j=1}^n
\left(
 \frac{1}{1+Y_j}
\right)^{c_{ij}}
=
\left(
\frac{Y_i}
{1+Y_i}\right)
=
\left(
\frac{Y_i}
{1+Y_i}\right)^2
\big(1+Y_i{}^{-1}\big),
\end{gather*}
which is equivalent to
\begin{gather}
\label{eq:Ynflat2}
Y_i^2=
\frac{\prod\limits_{j=1}^n (1+Y_j)^{2\delta_{ij}-c_{ij}}}
{
1+Y_i{}^{-1}
}.
\end{gather}

For simply laced types, $A_n$, $D_n$, $E_6$,
$E_7$, $E_8$, \eqref{eq:Ynflat2}  is obtained from
the {\em level} 3
constant $Y$-system~\eqref{eq:Y-ADE} of untwisted type $R_n^{(1)}$
by the specialization $Y^{(i)}_1=Y^{(i)}_2$ and identifying it with $Y_i$.
This is possible, due to the symmetry $Y^{(i)}_1\leftrightarrow
Y^{(i)}_2$ of level 3 $Y$-system~\eqref{eq:Y-ADE}.
(One can also view it as the folding of $A_2$ to $T_1$
in the second component of~$R_n\times A_2$.)
Since we discard the multiplicity in~\eqref{eq:CFdilog},
$L_{A^{\flat}}$~is the {\em half} of
 the right hand side of~\eqref{eq:dilogADE} with $\ell=3$.
 This agrees with Table~\ref{tab:L}.

Similarly,
for the rest of types,
\eqref{eq:Ynflat2}  is obtained from
the level~3 constant $Y$-system of type~$S_m^{(r)}$ or~$A_{2n}^{(2)}$
(in the sense of  \cite[Remark~9.22]{Inoue10c})
 by the specialization
 $Y^{(i)}_1=Y^{(i)}_2$,
and the latter is further obtained from
the level~3 constant $Y$-system of type $S_m^{(1)}$ or $A_{2n}^{(1)}$
 by the folding.
 Then, one can conf\/irm  Table~\ref{tab:L}.

Let us summarize the result.

\begin{theorem}
The identity \eqref{eq:CFdilog} holds.
Moreover,   except for type $T_n$,
the value $L_{A'}$ in~\eqref{eq:CFdilog}
has a~unified expression
\begin{gather*}
L_{A'}=\frac{mh^*}{h^*+\ell},
\end{gather*}
where $\ell=2$ for $A'=A$ and $\ell=3$ for $A'=A^{\flat}$,
and $m$ and $h^*$ are the rank  and Coxeter number
of~$S_m$ for the Langlands dual~$S_m^{(r)}$ of~$R_n^{(1)}$.
\end{theorem}

We remark that the dilogarithm identities for untwisted {\em and}
 nonsimply laced types
$B^{(1)}_n$, $C^{(1)}_n$, $F^{(1)}_4$, $G^{(1)}_2$ are also known~\cite{Inoue10a,Inoue10b}.
 It is natural to ask whether they will also appear from
Nil-DAHA.

\section{Connection to cluster algebraic method}

For the reader's convenience,
we include a brief explanation of the background
of the dilogarithm identity \eqref{eq:dilogADE},
 especially in  the cluster algebraic method.
 See \cite{Kashaev11} and references therein for more information.

{\em $(a)$ $Y$-systems and dilogarithm identities.}
As the name suggests, the constant $Y$-system
\eqref{eq:Y-ADE} is the constant version
of the following (non-constant) $Y$-system
\begin{gather}
\label{eq:Y-ADE2}
Y^{(i)}_m(u+1)
Y^{(i)}_m(u-1)
=
\frac{
\prod\limits_{j=1}^n \big(1+Y_m^{(j)}(u)\big)^{2\delta_{ij}-c_{ij}}}
{
\big(1+Y_{m-1}^{(i)}(u)^{-1}\big)
\big(1+Y_{m+1}^{(i)}(u)^{-1}\big)
},
\end{gather}
where the variables $Y^{(i)}_m(u)$ now carry the {\em spectral parameter}
$u\in \mathbb{C}$.
The $Y$-system~\eqref{eq:Y-ADE2} appears in
the thermodynamic Bethe ansatz
(TBA) analysis of the
deformation of conformal f\/ield theory.
A~constant solution $Y^{(i)}_m:=Y^{(i)}_m(u)$ of \eqref{eq:Y-ADE2},
which is constant with respect to~$u$,
satisf\/ies the constant $Y$-system~\eqref{eq:Y-ADE},
that also appears in the TBA analysis to calculate the ef\/fective
central charge of conformal f\/ield theory.

The $Y$-system \eqref{eq:Y-ADE2} has
the following two remarkable properties.

{\em $(i)$ Periodicity}
\begin{gather}
\label{eq:period1}
Y^{(i)}_m(u+2(h+\ell))=
Y^{(i)}_m(u).
\end{gather}

{\em $(ii)$ Dilogarithm identity.}
For any positive real solution of  \eqref{eq:Y-ADE2},
the following identity holds
\begin{gather}
\label{eq:dilogADE2}
\frac{6}{\pi^2}
\sum_{u=0}^{2(h+\ell)-1}
\sum_{i=1}^n \sum_{m=1}^{\ell-1}L\left(\frac{Y^{(i)}_m(u)}{1+Y^{(i)}_m(u)}\right) =
2(\ell-1)nh.
\end{gather}

The identity \eqref{eq:dilogADE}
is obtained by
specializing the identity
\eqref{eq:dilogADE2} to the (unique) positive real constant solution,
and dividing the both hand sides of \eqref{eq:dilogADE} by
the period $2(h+\ell)$.

{\em $(b)$ Cluster algebras and dilogarithm identities.}
It was once formidable to prove the properties~\eqref{eq:period1} and~\eqref{eq:dilogADE2} in full generality.
However, they are now proved and rather well understood
by the cluster algebraic method.
In general, to {\em any  period} of $y$-variables (coef\/f\/icients) of a cluster algebra,
the following dilogarithm identity is associated
\begin{gather}
\label{eq:dilogADE3}
\frac{6}{\pi^2}
 \sum_{t=1}^{p}L\left(\frac{y_{k_t}(t)}{1+y_{k_t}(t)}\right) =N_-.
\end{gather}
where $k_1, \dots, k_p$ are the sequence of mutations for which
 $y$-variables are periodic,
$y_{k_t}(t)$ is the  $y$-variable mutated at~$t$,
and~$N_-$ is the total number of $t\in \{ 1,\dots,p\}$ such that
the tropical sign of~$y_{k_t}(t)$ is minus.
One can apply this general result to our
$Y$-system~\eqref{eq:Y-ADE2}.
First, the $Y$-system is  embedded into the
$y$-variables of a certain cluster algebra.
Then, by proving the
periodicity  of $y$-variables of this cluster algebra,
we obtain the periodicity~\eqref{eq:period1} of the $Y$-system.
Finally, from the general identity~\eqref{eq:dilogADE3},
we obtain the dilogarithm identity~\eqref{eq:dilogADE2}
by calculating the constant term~$N_-$.

{\em $(c)$ Quantum cluster algebras and quantum dilogarithm identities.}
One can lift the result in~$(b)$ to the quantum case.
Namely, any period of $y$-variables of a cluster algebra can be
lifted to the period of {\em quantum $y$-variables} of the corresponding
{\em quantum cluster algebra}.
Then, to any such period, the quantum dilogarithm identity is
associated;
furthermore, taking the semiclassical limit of the quantum dilogarithm identity	
we recover the classical dilogarithm identity~\eqref{eq:dilogADE3}.
We expect that  the Rogers--Ramanujan type identities of~ \cite{Cherednik12}
 are also deduced from
 these quantum dilogarithm identities.
 Some result in this direction is obtained by~\cite{Cecotti10}.

\subsection*{Acknowledgements}
I thank Ivan Cherednik for raising this interesting question  to me.

\pdfbookmark[1]{References}{ref}
\LastPageEnding

\end{document}